\documentclass{IEEEtran4PSCC}
\usepackage[utf8]{inputenc}
\usepackage{graphicx}
\usepackage{multirow}
\usepackage{multicol}
\usepackage{lipsum}
\usepackage{amsmath}
\usepackage{amsbsy}
\usepackage{blindtext}
\usepackage{tcolorbox}
\usepackage{lscape}
\usepackage{amssymb}
\usepackage{scalerel}
\usepackage{mathabx}
\usepackage{verbatim}
\usepackage{booktabs}
\usepackage{rotating,tabularx}
\usepackage{mathrsfs} 
\usepackage[english]{babel}
\usepackage{url}
\usepackage{adjustbox}
\usepackage{soul}
\usepackage{xcolor}
\usepackage{cite}
\usepackage{tikz}
\usepackage{babel}
\usepackage[caption=false,font=footnotesize]{subfig}
\usepackage[T1]{fontenc}
\usepackage[utf8]{inputenc}
\usepackage{scalerel,stackengine}
\newcommand\reallywidecheck[1]{%
\savestack{\tmpbox}{\stretchto{%
  \scaleto{%
    \scalerel*[\widthof{\ensuremath{#1}}]{\kern-.6pt\bigwedge\kern-.6pt}%
    {\rule[-\textheight/2]{1ex}{\textheight}}
  }{\textheight}%
}{0.5ex}}%
\stackon[1pt]{#1}{\scalebox{-1}{\tmpbox}}%
}
\stackMath
\IEEEoverridecommandlockouts
\IEEEoverridecommandlockouts

\renewcommand{\baselinestretch}{0.981}
\usepackage[utf8]{inputenc}
\usepackage{array}
\usepackage{makecell}

\usepackage{bm}

\allowdisplaybreaks

\begin{document}

\title{\LARGE \bf Improving QC Relaxations of OPF Problems via Voltage Magnitude Difference Constraints and Envelopes for Trilinear Monomials}
\author{Mohammad Rasoul Narimani,$^{\ast}$ Daniel K. Molzahn,$^{\dagger}$  and Mariesa L. Crow$^{\ast}$%
\thanks{${\ast}$: ECE Dept., Missouri University of Science and Technology.}%
\thanks{${\dagger}$: Energy Systems Div., Argonne National Laboratory. Funding from the U.S. Department of
Energy, Office of Electricity Delivery and Energy Reliability under
award number \mbox{DE-AC-02-06CH11357}.}%
\thanks{The authors gratefully acknowledge insightful discussions with \mbox{Dr.~Carleton} \mbox{Coffrin} regarding the voltage magnitude difference constraints.}
}


\maketitle
\begin{abstract}
AC optimal power flow (AC~OPF) is a challenging non-convex optimization problem that plays a crucial role in power system operation and control. Recently developed convex relaxation techniques provide new insights regarding the global optimality of AC~OPF solutions. The quadratic convex (QC) relaxation is one promising approach that constructs convex envelopes around the trigonometric and product terms in the polar representation of the power flow equations. This paper proposes two methods for tightening the QC relaxation. The first method introduces new variables that represent the \emph{voltage magnitude differences} between connected buses. Using ``bound tightening'' techniques, the bounds on the voltage magnitude difference variables can be significantly smaller than the bounds on the voltage magnitudes themselves, so constraints based on voltage magnitude differences can tighten the relaxation. Second, rather than a potentially weaker ``nested McCormick'' formulation, this paper applies ``Meyer and Floudas'' envelopes that yield the convex hull of the trilinear monomials formed by the product of the voltage magnitudes and trignometric terms in the polar form of the power flow equations. Comparison to a state-of-the-art QC implementation demonstrates the advantages of these improvements via smaller optimality gaps.
\end{abstract}

 

\section{Introduction}
\label{introduction}
The optimal power flow (OPF) problem seeks an operating point that optimizes a specified objective subject to constraints from the network physics and engineering limits. Using the nonlinear AC power flow model to accurately represent the power flow physics results in the AC~OPF problem, which is non-convex, may have multiple local optima~\cite{bukhsh_tps}, and is generally NP-Hard~\cite{pascalNPhard,bienstock2015nphard}. 
A wide variety of algorithms have been applied in order to find locally optimal solutions~\cite{opf_litreview1993IandII,ferc4}.

Many recent research efforts have developed convex relaxations of OPF problems to obtain bounds on the optimal objective values, certify infeasibility, and, in some cases, achieve globally optimal solutions. Solutions from a relaxation are also useful for initializing certain local solution techniques~\cite{marley2016}. Convex relaxations are under active development with ongoing efforts aiming to improve the relaxations' computational tractability and tightness. Recent work is surveyed in~\cite{molzahn2017survey}.

The quadratic convex (QC) relaxation~\cite{coffrin2015qc} is one promising approach that uses convex envelopes around the trigonometric functions, squared terms, and bilinear products in the polar form of the power flow equations. The tightness of the QC relaxation depends on the size of the bounds on the voltage magnitude and angle difference variables. Therefore, bound tightening techniques, which use convex relaxations to infer tighter bounds than those initially specified in the OPF problem data, can improve the QC relaxation's tightness~\cite{coffrin2016strengthen_tps,chen2015,StrongSOCPRelaxations,arctan2}. Several enhancements have also been proposed to tighten the QC and other relaxations, including Lifted Nonlinear Cuts~\cite{coffrin2016strengthen_tps,chen2015cuts} that exploit voltage magnitude and angle difference bounds; tighter trigonometric envelopes~\cite{coffrin2016strengthen_tps,coffrin2016quadtrig} that leverage sign-definite angle difference bounds, which can sometimes be obtained via bound tightening; and a variety of valid inequalities, convex envelopes, and cutting planes~\cite{StrongSOCPRelaxations,arctan2}.


This paper proposes two additional improvements for tightening the QC relaxation. The first is based on the observation that adding redundant constraints to a non-convex optimization problem can tighten a relaxation~\cite{ruiz2011}. One approach for constructing appropriate constraints is to change coordinate systems. We derive constraints based on a coordinate change using voltage magnitude \emph{differences} in addition to the voltage magnitudes themselves. Bound tightening techniques are often more effective for variables representing voltage magnitude differences, thus resulting in tighter constraints.

The second improvement is related to the trilinear monomials formed by the product of the voltage magnitudes and the trignometric functions in the polar representation of the power flow equations. Previous formulations of the QC relaxation~\cite{coffrin2015qc,coffrin2016strengthen_tps} treat these monomials with recursive application of \mbox{McCormick} envelopes~\cite{mccormick1976}. While McCormick envelopes form the convex hull of \emph{bilinear} monomials, recursive application of McCormick envelopes does not necessarily yield the convex hulls of \emph{trilinear} monomials. We apply the potentially tighter envelopes developed by Meyer and Floudas~\cite{meyer2004a,meyer2004b}, which form the convex hulls of trilinear monomials.

This paper is organized as follows. Section~\ref{OPF_overview} overviews the OPF problem. Section~\ref{QC relaxation overview} reviews the QC relaxation of the OPF problem. Sections~\ref{voltagedifference} and~\ref{meyer&Floudasfacets} formulate our proposed improvements. Section~\ref{Numerical_results} evaluates the proposed improvements on various test cases. Section~\ref{conclusion} concludes the paper.

\section{Overview of Optimal Power Flow Problem}
\label{OPF_overview}
This section overviews the AC~OPF problem. Consider an $n$-bus system, where $\mathcal{N} = \left\lbrace 1, \ldots, n \right\rbrace$, $\mathcal{G}$, and $\mathcal{L}$ are the sets of buses, generators, and lines. Let ${P}_i^d + \mathbf{j}{Q}_i^d$ and ${P}_i^g + \mathbf{j} {Q}_i^g$ represent the active and reactive load demand and generation, respectively, at bus~$i\in\mathcal{N}$, where $\mathbf{j} = \sqrt{-1}$. Let $g_{sh,i} + \mathbf{j} b_{sh,i}$ denote the shunt admittance at bus $i$. Let $V_i$ and $\theta_i$ represent the voltage magnitude and angle at bus~$i\in\mathcal{N}$. For each generator $i\in\mathcal{G}$, define a quadratic generation cost function with coefficients $c_{2,i} \geq 0$, $c_{1,i}$, and $c_{0,i}$. Denote $\theta_{lm}=\theta_{l}-\theta_{m}$. Specified upper and lower limits are denoted by $\left(\overline{\,\cdot\,}\right)$ and $\left(\underline{\,\cdot\,}\right)$, respectively. Buses $i\in\mathcal{N}\setminus\mathcal{G}$ have generation limits set to zero. 

Each line $\left(l,m\right)\in\mathcal{L}$ is modeled as a $\Pi$ circuit with mutual admittance $g_{lm}+\mathbf{j} b_{lm}$ and shunt admittance $\mathbf{j} b_{sh,lm}$. (Our approach is applicable to more general line models, such the M{\sc atpower}~\cite{matpower} model that allows for off-nominal tap ratios and non-zero phase shifts.) Let ${P}_{lm}, {Q}_{lm}$, and $\overline{S}_{lm}$ represent the active and reactive power flows and the maximum apparent power flow limit on the line that connects buses~$l$ and~$m$.

Using these definitions, the OPF problem is
\begin{subequations}
\label{OPF formulation}
\begin{align}
\label{eq:objective}
& \min\quad \sum_{{i}\in \mathcal{G}}
c_{2i}\left(P_i^g\right)^2+c_{1i}\,P_i^g+c_{0i}\\ 
&\nonumber \text{subject to} \quad \left(\forall i\in\mathcal{N},\; \forall \left(l,m\right) \in\mathcal{L}\right) \\
\label{eq:pf1}
& P_i^g-P_i^d = g_{sh,i}\, V_i^2+\sum_{\substack{(l,m)\in \mathcal{L}\\ \text{s.t.} \hspace{3pt} l=i}} P_{lm}+\sum_{\substack{(l,m)\in \mathcal{L}\\ \text{s.t.} \hspace{3pt} m=i}} P_{ml}, \\
\label{eq:pf2}
& Q_i^g-Q_i^d = -b_{sh,i}\, V_i^2+\sum_{\substack{(l,m)\in \mathcal{L}\\ \text{s.t.} \hspace{3pt} l=i}} Q_{lm}+\sum_{\substack{(l,m)\in \mathcal{L}\\ \text{s.t.} \hspace{3pt} m=i}} Q_{ml},\\[-10pt]
\label{eq:OPF1}
& \theta_{ref}=0,\\
\label{eq:OPF2}
& \underline{P}_i^g\leq P_i^g\leq \overline{P}_i^g,\\
\label{eq:OPF3}
& \underline{Q}_i^g\leq Q_i^g\leq \overline{Q}_i^g,\\
\label{eq:OPF4}
& \underline{V}_i\leq V_{i} \leq \overline{V}_i,\\
\label{eq:OPF5}
&\underline{\theta}_{lm}\leq \theta_{lm}\leq \overline{\theta}_{lm},\\
\label{eq:pik}
&P_{lm} = g_{lm} v_l^2 - g_{lm} V_l V_m\cos\left(\theta_{lm}\right) - b_{lm} V_l V_m\sin\left(\theta_{lm}\right), \\[-8pt]
\label{eq:qik}
& \nonumber Q_{lm} = -\left(b_{lm}+b_{sh,lm}/2\right) V_l^2 + b_{lm} V_l V_m\cos\left(\theta_{lm}\right)\\ &\qquad\qquad  - g_{lm} V_l V_m\sin\left(\theta_{lm}\right), \\
\label{eq:OPF8}
& \left(P_{lm}\right)^2+\left(Q_{lm}\right)^2 \leq \left(\overline{S}_{lm}\right)^2, \\
\label{eq:OPF9}
& \left(P_{ml}\right)^2+\left(Q_{ml}\right)^2 \leq \left(\overline{S}_{lm}\right)^2.
\end{align}
\end{subequations}
The objective function~\eqref{eq:objective} minimizes the active power generation cost. Constraints~\eqref{eq:pf1} and~\eqref{eq:pf2} enforce power balance at each bus. Constraint~\eqref{eq:OPF1} sets the angle reference. Constraints~\eqref{eq:OPF2}--\eqref{eq:OPF5} limit the active and reactive power generation, voltage magnitudes, and angle differences between connected buses. Constraints~\eqref{eq:pik}--\eqref{eq:qik} relate the voltage phasors and power flows on each line, and \eqref{eq:OPF8}--\eqref{eq:OPF9} limit the apparent power flows into both terminals of each line.
 
\section{Review of the QC Relaxation}
\label{QC relaxation overview}

\subsection{Formulation of the QC Relaxation}
The QC relaxation is formed by defining new variables $w_{ii}$, $w_{lm}$, $c_{lm}$, and $s_{lm}$ for the products of voltage magnitudes and the trilinear monomials representing the products of voltage magnitudes and trignometric functions for connected buses:
\begin{subequations}
\label{eq:cs}
\begin{align}
w_{ii} &= V_i^2,  & \forall i \in\mathcal{N}, \\
w_{lm} &= V_l V_m,  & \forall \left(l,m\right) \in\mathcal{L}, \\
c_{lm} & =  w_{lm} \cos\left(\theta_{lm} \right), & \forall \left(l,m\right) \in\mathcal{L}, \\
s_{lm} & = w_{lm}\sin\left(\theta_{lm} \right),  & \forall \left(l,m\right) \in\mathcal{L}.
\end{align}
\end{subequations}
For each line $\left(l,m\right)\in\mathcal{L}$, these definitions imply the following relationships between the variables $w_{ll}$, $c_{lm}$, and $s_{lm}$:
\begin{subequations}
\label{eq:cs_relationships}
\begin{align}
\label{eq:Jabr}
&c_{lm}^2+s_{lm}^2=w_{ll}w_{mm},\\
\label{eq:cs_relationships_c}
&c_{lm}=c_{ml}, \\
\label{eq:cs_relationships_s}
&s_{lm}=-s_{ml}
\end{align}
\end{subequations}

The QC relaxation is formulated by enclosing the squared and bilinear product terms in convex envelopes, here represented as set-valued functions:
\begin{subequations}
\label{eq:product_envelopes}
\begin{align}
\label{eq:squareenvelopes}
\langle x^2\rangle^T =
\begin{cases}
\widecheck{x}: \begin{cases}\check{x} \geq x^2,\\
\widecheck{x} \leq \left({\overline{x}+\underline{x}}\right) x-{\overline{x} \underline{x}}.\\
\end{cases}
\end{cases}\\
\label{eq:mccormick}
\langle {xy}\rangle^M  =
\begin{cases}
\widecheck{xy}:\begin{cases}
\widecheck{xy} \geq {\underline{x}} y+ {\underline{y}} x-{\underline{x} \underline{y}},\\
\widecheck{xy} \geq {\overline{x}} y+ {\overline{y}} x-{\overline{x} \overline{y}},\\
\widecheck{xy} \leq {\underline{x}} y+ {\overline{y}} x-{\underline{x}} {\overline{y}},\\
\widecheck{xy} \leq {\overline{x}} y+ {\underline{y}} x-{\overline{x} \underline{y}}.\\
\end{cases}
\end{cases}
\end{align}
\end{subequations}
where $\widecheck{x}$ and $\widecheck{xy}$ are ``dummy'' variables representing the corresponding set.
The envelope $\langle x^2\rangle^T$ is the convex hull of the square function. The so-called ``McCormick envelope'' $\langle {xy}\rangle^M$ is the convex hull of a bilinear product~\cite{mccormick1976}.

The QC relaxation also formulates convex envelopes $\left\langle\sin\left(x\right) \right\rangle^S$ and $\left\langle\cos\left(x\right) \right\rangle^C$ for the trigonometric functions:
\begin{subequations}
\label{eq:convex_envelopes_sin&cos}
\begin{align}
\label{eq:sine envelope}
\nonumber &\left\langle \sin(x)\right\rangle^S =\\
&\quad\;\begin{cases}
\widecheck{S}:\begin{cases}
\widecheck{S}\leq\cos\left(\frac{x^m}{2}\right)\left(x-\frac{x^m}{2}\right)+\sin \left(\frac{x^m}{2}\right),\\
\widecheck{S}\geq\cos\left(\frac{x^m}{2}\right)\left(x+\frac{x^m}{2}\right)-\sin\left(\frac{x^m}{2}\right),\\
\widecheck{S}\geq\frac{\sin\left({\underline{x}}\right)-\sin\left(\overline{x}\right)}{{\underline{x}-\overline{x}}}\left(x-{\underline{x}}\right)+\sin\left({\underline{x}}\right) \text{if~} \underline{x}\geq0,\\
\widecheck{S}\leq\frac{\sin\left({\underline{x}}\right)-\sin\left({\overline{x}}\right)}{{\underline{x}-\overline{x}}}\left(x-{\underline{x}}\right)+\sin\left({\underline{x}}\right) \text{if~} {\overline{x}}\leq0.
\end{cases}
\end{cases}\\
\label{eq:cosine envelope}
&\nonumber\left\langle\cos(x)\right\rangle^C=\\
&\quad\;\begin{cases}
\widecheck{C}:\begin{cases}
\widecheck{C}\leq 1-\frac{1-\cos\left({x^m}\right)}{\left(x^m\right)^2}x^2,\\
\widecheck{C}\geq\frac{\cos\left(\underline{x}\right)-\cos\left({\overline{x}}\right)}{{\underline{x}-\overline{x}}}\left(x-{\underline{x}}\right)+\cos\left({\underline{x}}\right).  
\end{cases}
\end{cases}
\end{align}
\end{subequations}
where $x^m= \max(\left|\underline{x}\right|,\left|\overline{x}\right|)$. The dummy variables $\check{S}$ and $\check{C}$ again represent the corresponding set. For $-90^\circ < \underline{x} < \overline{x} < 90^\circ$, bounds on the sine and cosine functions are
%
\begin{subequations}
\begin{align}
& \underline{s} = \sin\left(\underline{x}\right) \leq \sin(x) \leq \overline{s} = \sin\left(\overline{x}\right), \\
\nonumber & \underline{c} = \min\left(\cos(\underline{x}),\cos(\overline{x})\right) \leq \cos(x) \\ & \quad \leq \overline{c} \!=\! \begin{cases} \max\left(\cos(\underline{x}),\cos(\overline{x})\right),\; \text{if~} \mathrm{sign}\left(\underline{x}\right) \!=\!  \mathrm{sign}\left(\overline{x}\right), \\ 1, \text{~otherwise}. \end{cases}\raisetag{1em}
\end{align}
\end{subequations} 

Slightly abusing notation, the QC relaxation is formed by replacing the square, product, and trigonometric terms in~\eqref{OPF formulation} with the variables $w_{ii}$, $w_{lm}$, $c_{lm}$, and $s_{lm}$ in these envelopes:
\begin{subequations}
\label{eq:qc}
\begin{align}
& \min \quad \sum_{{i}\in \mathcal{G}} c_{2i}\left(P_i^g\right)^2+c_{1i}\,P_i^g+c_{0i}\\ 
&\nonumber \text{subject to} \quad \left(\forall i\in\mathcal{N},\; \forall \left(l,m\right) \in\mathcal{L}\right) \\
\label{eq:qc_p}
& P_i^g-P_i^d = g_{sh,i}\, w_{ii}+\sum_{\substack{(l,m)\in \mathcal{L}\\ \text{s.t.} \hspace{3pt} l=i}} P_{lm}+\sum_{\substack{(l,m)\in \mathcal{L}\\ \text{s.t.} \hspace{3pt} m=i}} P_{ml}, \\
\label{eq:qc_q}
& Q_i^g-Q_i^d = -b_{sh,i}\, w_{ii}+\sum_{\substack{(l,m)\in \mathcal{L}\\ \text{s.t.} \hspace{3pt} l=i}} Q_{lm}+\sum_{\substack{(l,m)\in \mathcal{L}\\ \text{s.t.} \hspace{3pt} m=i}} Q_{ml},\\
\label{eq:qc_V}
&  (\underline{V}_i)^2\leq w_{ii} \leq (\overline{V}_i)^2,\\
\label{eq:qc_pik}
& P_{lm} = g_{lm} w_{ll} - g_{lm} c_{lm} - b_{lm} s_{lm}, \\
\label{eq:qc_qik}
& Q_{lm} = -\left(b_{lm}+b_{sh,lm}/2\right) w_{ii} + b_{lm} c_{lm}- g_{lm} s_{lm}, \\
\label{eq:qc_wii}
&  w_{ii} \in\left\langle V_i^2 \right\rangle^T, \\
\label{eq:qc_wik}
&  w_{lm} \in \left\langle V_l V_m \right\rangle^M, \\
\label{eq:qc_cik}
&  c_{lm} \in \left\langle w_{lm}\left\langle\cos\left(\theta_{lm} \right)\right\rangle^C\right\rangle^M, \\
\label{eq:qc_sik}
& s_{lm} \in \left\langle w_{lm} \left\langle\sin\left(\theta_{lm} \right)\right\rangle^S\right\rangle^M, \\
\label{eq:qc_jabr}
&  c_{lm}^2+s_{lm}^2 \leq w_{ll}\,w_{mm} \\
\label{eq:qc_others}
&  \text{Equations~}\eqref{eq:OPF1}\text{--}\eqref{eq:OPF5},\,\eqref{eq:OPF8}\text{--}\eqref{eq:OPF9},\,\eqref{eq:cs_relationships_c},\,\eqref{eq:cs_relationships_s}.
\end{align}
\end{subequations}
Note that the non-convex constraint~\eqref{eq:Jabr} is relaxed to~\eqref{eq:qc_jabr} using a less-stringent rotated second-order cone constraint~\cite{Jabr2006}. Also note that the trilinear terms in~\eqref{eq:pik} and~\eqref{eq:qik} are addressed in~\eqref{eq:qc_wik}--\eqref{eq:qc_sik} by recursively applying \mbox{McCormick} envelopes~\eqref{eq:mccormick} (i.e., first applying~\eqref{eq:mccormick} to the product of voltage magnitudes to obtain $w_{lm}$ and then to the product of $w_{lm}$ and $\left\langle\cos\left(\theta_{lm} \right)\right\rangle^C$ or $\left\langle\sin\left(\theta_{lm} \right)\right\rangle^S$). 

The optimization problem~\eqref{eq:qc} is a second-order cone program (SOCP), which is convex and can be solved efficiently using commercial tools (e.g., CPLEX, Gurobi, and Mosek).

\subsection{Bound Tightening and Other Improvements}
\label{subsec:bt}
The tightness of the QC relaxation strongly depends on the accuracy of the bounds on voltage magnitudes, $\underline{V}_i,~ \overline{V}_i$, and angle differences, $\underline{\theta}_{lm},~\overline{\theta}_{lm}$. The values specified in the dataset for these bounds may be significantly larger than the values that are actually achievable due to the restrictions imposed by other constraints. In other words, certain bounds may never be binding. Exploiting this observation, bound tightening algorithms yield tighter bounds that improve the QC relaxation~\cite{coffrin2016strengthen_tps,StrongSOCPRelaxations,chen2015,arctan2}.

We apply the optimization-based bound tightening algorithm in~\cite{coffrin2016strengthen_tps}, which iteratively minimizes and maximizes each voltage magnitude and angle difference variable subject to the QC relaxation's constraints. For instance, consider the upper bound on the voltage magnitude at bus~$1$:
\begin{align}
\label{eq:bt_example}
w_{11}^\ast = \max \quad w_{11} \quad \text{subject to}\quad \eqref{eq:qc_p}\text{--}\eqref{eq:qc_others}.
\end{align}
The value $w_{11}^\ast$ upper bounds the maximum achievable value of $\left(V_{1}\right)^2$ within the feasible space. If $w_{11}^\ast < \left(\overline{V}_{1}\right)^2$, then~\eqref{eq:bt_example} provides a smaller value of $\sqrt{w_{11}^\ast}$ for the upper bound on $V_1$, which tightens the QC relaxation. Since tightening the bound on any variable may improve the achievable bounds on other variables, the bound tightening algorithm proceeds iteratively until no further bounds can be tightened. Optimization-based bound tightening algorithms, e.g.,~\cite{StrongSOCPRelaxations,coffrin2016strengthen_tps,arctan2}, are typically slower than analytical methods~\cite{chen2015} but provide tighter bounds.

Previous literature proposes a variety of other improvements to the QC relaxation. To form a benchmark for comparing our improvements, we augment~\eqref{eq:qc} with quadratic envelopes for the trigonometric terms~\cite{coffrin2016quadtrig}, arctangent envelopes~\cite{arctan2}, and lifted nonlinear cuts (LNC)~\cite{coffrin2016strengthen_tps,chen2015cuts}.

\section{Voltage Magnitude Difference Constraints}
\label{voltagedifference}



As discussed in Section~\ref{subsec:bt}, the QC relaxation's tightness strongly depends on having accurate bounds on voltage magnitudes and angle differences. While bound tightening techniques are often successful in reducing the range of the phase angle differences, tightening the voltage magnitudes can be more challenging since OPF feasible spaces typically contain points for which the voltage magnitudes are both near the top and near the bottom of their allowed ranges. The bound tightening algorithms are therefore often unable to significantly improve the voltage magnitude bounds. 

However, there is usually an exploitable correlation between the voltage magnitudes at neighboring buses. While the voltage magnitudes at a pair of neighboring buses may be near their upper limits or near their lower limits, typical problems with limited reactive power injection capabilities require that these voltage magnitudes must be \emph{close to each other}. This suggests that ``box constraints'' on the voltage magnitudes~\eqref{eq:OPF4} are not a good match to the voltage magnitude variation exhibited in typical OPF feasible spaces. 

\begin{figure}[!t]
	\centering
	\includegraphics[scale=0.36,,trim={0cm 0.75cm 0cm 0},clip]{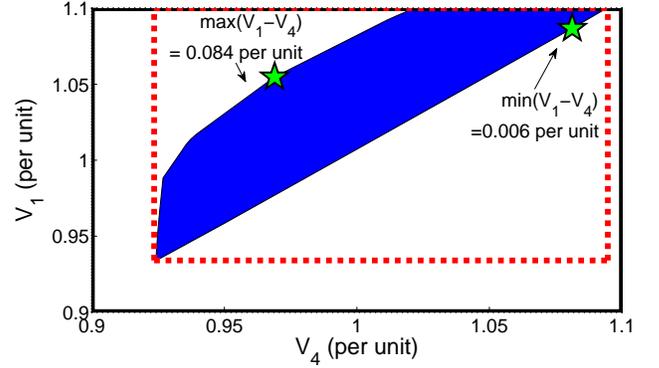}
	\caption{A projection of the feasible space for the ``case6\_c''~\cite{nesta3} test system.}
	\label{fig_sim}
	\vspace*{-1.3em}
\end{figure}

As an illustrative example, Figure~\ref{fig_sim} shows a projection of the feasible space, generated using the approach in~\cite{molzahn-opf_spaces}, for the six-bus system ``case6\_c''~\cite{nesta3} in terms of certain voltage magnitudes. The ranges of the voltage magnitude variations after implementing a bound tightening approach are shown by the dashed lines. The best achievable voltage magnitude bounds are only 17.0\% tighter than the originally specified bounds for this case. In contrast, as shown in Figure~\ref{fig_sim}, the \emph{difference} in voltage magnitudes between neighboring buses can be significantly tighter (e.g., 80.5\% tighter for the example in Figure~\ref{fig_sim}). 

To exploit this observation, we derive new constraints by representing the decision variables in an alternate coordinate system. Let $\mathbf{A}_{inc} \in\mathbb{R}^{\left|\mathcal{L}\right|\times \left|\mathcal{N}\right|}$ denote the network incidence matrix, which has rows corresponding to the lines $\left(l,m\right)\in\mathcal{L}$ with $+1$ in the $i^{th}$ entry and $\text{-}1$ in the $k^{th}$ entry. Define $V^\Delta\in\mathbb{R}^{\left|\mathcal{L}\right|}$ as the vector of voltage differences between neighboring buses, $V^\Delta = \mathbf{A}_{inc} V$ (i.e., $V^\Delta_{lm} = V_l - V_m$). Rewriting the voltage magnitude products $V_l V_m$ using $V^\Delta$ yields
\begin{align}
\label{eq:qc_aug}
V_l\, V_m = \left(V_l^2 + V_m^2 - \left(V^\Delta_{lm}\right)^2\right)/\,2.
\end{align}
Applying the envelopes in~\eqref{eq:squareenvelopes} for each term in~\eqref{eq:qc_aug} gives
\begin{subequations}
\label{eq:qc_aug_relax}
\begin{align}
w_{lm} & = \left(w_{ll} + w_{mm} - W^\Delta_{lm}\right)/\,2, \\
W^\Delta_{lm} & \in \left\langle \left(V^\Delta_{lm}\right)^2 \right\rangle^T.
\end{align}
\end{subequations}

A valid inequality is also formed by expanding $\left(V_{l} - V_{m}\right)^2$:
\begin{align}
\label{eq:qc_aug_inequ}
\left(V^\Delta_{lm}\right)^2 \leq V_l^2 - 2\, V_l\, V_m + V_m^2.
\end{align}
Relaxing~\eqref{eq:qc_aug_inequ} using~\eqref{eq:product_envelopes} yields
\begin{align}
\label{eq:qc_aug_inequ_relax}
\left(V^\Delta_{lm}\right)^2 \leq w_{ll} - 2w_{lm}  + w_{mm}.
\end{align}
Note that it is not necessary to use a convex envelope on the term $V^\Delta_{lm}$ since~\eqref{eq:qc_aug_inequ_relax} is already an SOCP constraint.

Finally, we leverage the relaxation proposed in~\cite{taylorlinearRelaxationOPF}, which is derived by taking linear combinations of the non-linear expressions for the active and reactive line flow expressions. Specifically, the following constraint from~\cite{taylorlinearRelaxationOPF} couples the voltage magnitude differences and the power flows:
\begingroup
\allowdisplaybreaks[0]
\begin{align}
\label{eq:consistency}
 V_l^2 - V_m^2 = \left(\frac{g_{lm}\left(P_{lm} - P_{ml}\right) - b_{lm}\left(Q_{lm} - Q_{ml}\right)}{g_{lm}^2 + b_{lm}^2 + b_{lm}\frac{b_{sh,lm}}{2}}\right).
\end{align}
\endgroup

Factoring the left hand side of~\eqref{eq:consistency} yields $V_l^2 - V_m^2 = V^\Delta_{lm} \left(V_l + V_m\right)$. Relaxing this expression yields
\begin{subequations}
\label{eq:consistency_relax}
\begingroup
\allowdisplaybreaks[0]
\begin{align}
\nonumber
&w_{ll} - w_{mm} = \hat{W}_{lm,l} + \hat{W}_{lm,m} \label{eq:consistency_relax_main}  \\ & \qquad = \left(\frac{g_{lm}\left(P_{lm} - P_{ml}\right) - b_{lm}\left(Q_{lm} - Q_{ml}\right)}{g_{lm}^2 + b_{lm}^2 + b_{lm}\frac{b_{sh,lm}}{2}}\right), \\
& \hat{W}_{lm,l} \in \left\langle V^\Delta_{lm} \, V_l \right\rangle^M,\\
& \hat{W}_{lm,m} \in \left\langle V^\Delta_{lm} \, V_m \right\rangle^M.
\end{align}
\endgroup
\end{subequations}
Observe that~\eqref{eq:consistency_relax_main} describes two constraints.

Our proposed improvement based on voltage magnitude differences augments the QC relaxation~\eqref{eq:qc} with constraints~\eqref{eq:qc_aug_relax}, \eqref{eq:qc_aug_inequ_relax}, and \eqref{eq:consistency_relax}.
The main advantage of these constraints is the quality of the achievable bounds on the voltage magnitude differences $V_{lm}^{\Delta}$. These bounds are computed by extending the bound tightening techniques described in Section~\ref{subsec:bt} to directly consider to the variables $V_{lm}^{\Delta}$. This requires initially specified bounds on $V_{lm}^{\Delta}$, which are derived from the bounds on the voltage magnitudes, $\underline{V}_l$, $\overline{V}_l$, $\underline{V}_m$, and $\overline{V}_m$:
\begin{equation}
\label{eq:initial_bounds}
\underline{V}_l - \overline{V}_m \leq V^\Delta_{lm} \leq \overline{V}_l - \underline{V}_m.
\end{equation}
After applying bound tightening to the voltage magnitudes, voltage angle differences, and voltage magnitude differences, bounds on the remaining variables ($W_{lm}$, $\hat{W}_{lm,l}$, and $\hat{W}_{lm,m}$) are derived by straightforward manipulations of the bounds on voltage magnitudes and voltage magnitude differences.

\section{Trilinear Envelopes}
\label{meyer&Floudasfacets}
Previous formulations of the QC relaxation recursively apply McCormick envelopes~\eqref{eq:mccormick} to represent the trilinear products formed by the voltage magnitudes and trigonometric terms. However, this approach rarely results in the convex hull of the trilinear products~\cite{meyer2004a}. The Meyer and Floudas envelopes~\cite{meyer2004a,meyer2004b} form the convex hulls of trilinear products whose variables range in a box. These envelopes thus provide a mechanism for strengthening the QC relaxation.

Due to the signs of the variables (i.e., positive voltage magnitudes and cosine terms, sign-indefinite sine terms), only certain facets of these envelopes are applicable to the QC relaxation. 
The appendix provides equations for these facets.


\section{Numerical Results}
\label{Numerical_results}

This section demonstrates the proposed improvements using test cases from the NESTA 0.7.0 archive~\cite{nesta3} and four cases ``nmwc14'', ``nmwc24,'' ``nmwc57,'' and ``nmwc118'' from~\cite{ACCpaper}. With large optimality gaps between the objective values from the best known local optima and the lower bounds from various relaxations, these test cases challenge a variety of solution algorithms and are therefore suitable for our purposes.


The implementation uses MATLAB~2013a, YALMIP \mbox{2016.09.30}~\cite{yalmip}, Mosek \mbox{8.0.0.42}, and a laptop computer with an i5 3.20~GHz processor and 8~GB of RAM.

Table~\ref{tab:QC_results} details the results for selected test cases.
The first column indicates the test case. The second column provides the objective value from M{\sc atpower}~\cite{matpower}.
The next group of columns presents the optimality gaps corresponding to the solution of a QC relaxation variant relative to the local solution from M{\sc atpower}. 
The optimality gap~is
\begin{align}
\label{eq:optimality gap}
\text{\emph{Optimality~gap}}=\left(\dfrac{\text{\emph{Local~solution}} - \text{\emph{QC bound}}}{\text{\emph{QC bound}}}\right).
\end{align}
For many applications, such as branching algorithms that compute global optima~\cite{StrongSOCPRelaxations,arctan2,chen2015cuts}, mixed-integer problems~\cite{coffrin2016quadtrig,kocuk2017_ots}, and certain bi-level problems~\cite{molzahn2018pscc}, the optimality gap is of primary importance. We therefore use the optimality gap to measure the relaxations' tightness.\footnote{Note that the optimality gap depends on both the lower bound from the relaxation and the upper bound from a local solution. Thus, non-zero gaps may be partially due to a suboptimal local solution. However, the same local optima are used to compute the optimality gap for each relaxation, and the gaps can therefore be consistently compared among various relaxations for each test case.}
The final group of columns in Table~\ref{tab:QC_results} provides the solution times, listing both the bound tightening time and the QC relaxation's execution time. Note that the bounds were tightened using the corresponding variant of the QC relaxation in the computations. For typographical purposes, Table~\ref{tab:QC_results} uses several abbreviations: ``All Constraints'' (\emph{All Cons.}), ``without'' (\emph{w/o}), ``Meyer and Floudas Envelopes'' (\emph{MF}), ``Voltage Magnitude Difference constraints'' ($\Delta$), and ``Bound Tightening'' (\emph{BT}).


The results indicate that bound tightening has a substantial impact on the optimality gaps for all variants of the QC relaxation. For instance, comparing the third and seventh columns in Table~\ref{tab:QC_results} reveals that applying bound tightening reduces the gaps for ``nesta$\_$case30$\_$fsr$\_\_$api'' and ``nesta$\_$case118$\_$ieee$\_\_$api'' cases by $78.45$\% and $63.45$\%, respectively. The reinforces the fact that the accuracy of the QC relaxation strongly depends on the tightness of the bounds. 

Comparing the fourth and seventh columns with the third column demonstrates the impact of the Meyer and Floudas envelopes and voltage difference constraints, both individually and jointly. For instance,  
the optimality gap for ``nesta$\_$case118$\_$ieee$\_\_$api'' without applying these constraints was $18.71$\% while applying the Meyer and Floudas envelopes and the voltage difference constraints reduces the gap to $13.9$\% and $18.32$\%, respectively. Applying both at the same time reduces the gap to $13.61$\%, revealing that the Meyer and Floudas envelopes are the larger contributor to the improvement for this test case.  
Similar results are obtained for ``nesta$\_$case30$\_$fsr$\_\_$api''. Without the Meyer and Floudas envelopes and the voltage difference constraints, the gap is $5.50$\%. Applying these improvements reduces the gap by $1.34$\% and $0.21$\%, respectively. For most of the case studies in Table~\ref{tab:QC_results}, the Meyer and Floudas envelopes are responsible for more of the improvement than the voltage difference constraints. However, there are cases where the opposite is true, such as ``nmwc118'', ``nmwc57'', and ``nmwc14''. For these cases, the voltage difference constraints outperformed the Meyer and Floudas envelopes in reducing the optimality gap, by up to $5.89$\% in the case of ``nmwc118''.

The results suggest that the Meyer and Floudas envelopes and the voltage difference constraints are most effective when applied in combination with a bound tightening algorithm. However, there are cases, such as ``nesta$\_$case73$\_$ieee$\_$rts$\_\_$api'' and ``nesta$\_$case29$\_$edin$\_\_$sad'' where the proposed improvements have significant impact even without bound tightening ($6.39$\% and $8.30$\% reductions, respectively). Note that the Meyer and Floudas envelopes play a more important role in both cases. For instance, they reduce the optimality gap for ``nesta$\_$case29$\_$edin$\_\_$sad'' by almost $8.30$\%,  whereas the voltage difference constraints only reduce the gap by $0.01$\%. This matches the intuition that the voltage magnitude difference constraints strongly depend on tight bounds on $V_{lm}^\Delta$.

Several comparisons underscore the contributions of different improvements to a basic QC relaxation (with no previous or proposed improvements, i.e., without applying bound tightening, the approaches proposed in this paper, or those in~\cite{coffrin2016quadtrig,StrongSOCPRelaxations,arctan2,coffrin2016strengthen_tps,chen2015cuts}). 
Separately adding different improvements to the basic QC relaxation reveals the individual contributions. The optimality gap of the basic QC relaxation for ``nesta$\_$case73$\_$ieee$\_$rts$\_\_$api'' is $16.50$\%. Separately adding the LNC constraints in~\cite{coffrin2016strengthen_tps,chen2015cuts} and the arctangent envelopes in~\cite{arctan2}  
does not reduce the gap while separately adding the voltage difference constraints and the Meyer and Floudas envelopes reduces the gap by $0.02$\%, and $5.87$\%, respectively. Note that using bound tightening with the basic QC relaxation reduces the gap by $9.05$\%. Similarly, the optimality gap resulting from applying the basic QC relaxation to ``nesta$\_$case29$\_$edin$\_\_$sad'' is $34.53$\%. Separately enforcing the LNC constraints and the voltage magnitude difference constraints does not reduce the gap while the arctangent envelopes and the Meyer and Floudas envelopes reduce the gap by $6.58$\% and $14.89$\%, respectively. For this case, it is interesting to note that the bound tightening approach alone only reduces the gap by $0.61$\%.

The impact of the voltage magnitude difference constraints strongly depends the quality of the bounds on $V_{lm}^\Delta$. Thus, applying these constraints without using bound tightening has a limited effect, as discussed above. In contrast, the voltage magnitude difference constraints contribute to reducing the optimality gap when combined with a bound tightening approach. For instance, these constraints reduce the optimality gap for ``nmwc118'' by $6.08$\%, whereas the Meyer and Floudas envelopes only reduce the gap by $0.19$\%. 
Thus, the contributions of each improvement to reducing the optimality gap depend on the test case. Our future work includes identifying which system characteristics are most relevant for various types of improvements. 

Our proposed improvements substantially reduce the optimality gaps for many challenging test cases. As shown in Table~\ref{tab:QC_results}, this improved tightness comes at the cost of slower (but still tractable) computational times for some test cases. 
Comparing the last two columns in Table~\ref{tab:QC_results} reveals that enforcing the Meyer and Floudas envelopes and the voltage difference constraints results in less than a $43.6$\% increase in the time required to solve the QC relaxation (without bound tightening) on average across the test cases. 
Comparing the execution times in the ninth and twelfth columns of Table~\ref{tab:QC_results} shows that adding the Meyer and Floudas envelopes and the voltage difference constraints has a disparate impact on the total execution time (bound tightening plus QC execution). There are cases such as ``nesta$\_$case29$\_$edin$\_\_$sad'' where enforcing these constraints reduces the execution time by $64.9$\%. For these cases, the bound tightening algorithm converges in fewer iterations, which more than accounts for the additional time required per iteration due to the addition of new variables and constraints. Since the bound tightening times dominate the execution time for the QC relaxation, the overall time decreases for some cases. Conversely, other test cases require more time, resulting in an average increase of $14.6$\% over all the test cases and up to an $75.0$\% increase for some cases.


\section{Conclusion}
\label{conclusion}
This paper proposes and empirically tests two improvements for the QC relaxation of the OPF problem: a set of constraints based on voltage magnitude differences and the Meyer and Floudas envelopes for trilinear monomials. The former relies on the observation that bound tightening algorithms can effectively tighten the voltage magnitude differences between connected buses. The latter yields the convex hull of the trilinear monomials in contrast to the potentially weaker nested McCormick formulation used in previous work. Comparison to a state-of-the-art QC implementation demonstrates the value of these improvements via reduced optimality gaps on challenging test cases while maintaining computational tractability. Our ongoing work aims to improve computational speed by targeting the application of the bound tightening techniques to the most relevant variables. Other ongoing work is developing further improvements to convex relaxations based on physically intuitive coordinate transformations.


\appendix[Expressions for the Meyer and Floudas Envelopes]
\label{appendix:Meyer&Floudas}
This appendix provides the facets of the Meyer and Floudas envelopes that are applicable to the QC relaxation~\eqref{eq:qc}. In the following seven boxes, the upper portion gives conditions for which the constraints in the lower portion apply. 

We define $\widecheck{S} \in \left\langle \sin \left(\theta_{lm}\right) \right\rangle^S$, where this trigonometric envelope is given in~\eqref{eq:sine envelope}, and $V_i$ as the voltage magnitude at bus~$i$ as in~\eqref{eq:qc}. 
Let $\left\langle x\,y\,z  \right\rangle^{MF}$ denote the convex hull defined by the Meyer and Floudas envelopes for the trilinear product of three generic variables, $x$, $y$, and $z$. The variable $\widecheck{s}_{lm}\in \left\langle V_l\, V_m\, \widecheck{S} \right\rangle^{MF}$ replaces $s_{lm}$ in~\eqref{eq:qc}. Note that multiple cases may apply simultaneously (e.g., Case~IV implies Case~I). The same procedure is applied using $\widecheck{C} \in \left\langle \cos \left(\theta_{lm}\right) \right\rangle^C$, with the variable  $\widecheck{c}_{lm}\in \left\langle V_l\, V_m\, \widecheck{C} \right\rangle^{MF}$ replacing $c_{lm}$ in~\eqref{eq:qc}. Since the cosine function is non-negative in the first and fourth quadrants, only Cases~II and~III are applicable for this function.

\begin{figure}
\subfloat{%
\begin{tcolorbox}[width=\columnwidth,colback={white},title={\footnotesize Case~I:~$\overline{s}\leq 0$.},colbacktitle=white,coltitle=black,boxsep=0.2mm]   
\footnotesize
   \begin{subequations}
\begin{align*}
\\[-25pt]\widecheck{s}_{lm}&\geq \overline{V}_m \underline{s} V_l+ \underline{V}_l \underline{s} V_m+ \underline{V}_l \underline{V}_m\widecheck{S}-\underline{V}_l \overline{V}_m \underline{s}- \underline{V}_l \underline{V}_m \underline{s},\\ 
\widecheck{s}_{lm}&\geq \overline{V}_m \underline{s} V_l+\underline{V}_l \overline{s} V_m+ \underline{V}_l \overline{V}_m\widecheck{S}-\underline{V}_l \overline{V}_m \underline{s}- \underline{V}_l \overline{V}_m \overline{s},\\
\widecheck{s}_{lm}&\geq \underline{V}_m \overline{s} V_l+\overline{V}_l \underline{s} V_m+\overline{V}_l \underline{V}_m\widecheck{S}-\overline{V}_l \underline{V}_m \overline{s}- \overline{V}_l \underline{V}_m \underline{s},\\
\widecheck{s}_{lm}&\geq \underline{V}_m \overline{s} V_l+\overline{V}_l \overline{s} V_m+\overline{V}_l \overline{V}_m\widecheck{S}-\overline{V}_l \underline{V}_m \overline{s}- \overline{V}_l \overline{V}_m \overline{s},\\
\widecheck{s}_{lm}&\geq \underline{V}_m \underline{s} V_l+\overline{V}_l \underline{s} V_m+ \underline{V}_l \underline{V}_m\widecheck{S}-\overline{V}_l \underline{V}_m \underline{s}- \underline{V}_l \underline{V}_m \underline{s},\\
\widecheck{s}_{lm}&\geq \overline{V}_m \overline{s} V_l+\underline{V}_l \overline{s} V_m+\overline{V}_l \overline{V}_m\widecheck{S}-\overline{V}_l \overline{V}_m \overline{s}- \underline{V}_l \overline{V}_m \overline{s}.
\end{align*}
\end{subequations}
\end{tcolorbox}
}\\[0.75em]
\subfloat{%
\begin{tcolorbox}[width=\columnwidth,colback={white},title={\footnotesize Case~II: ~$\underline{s}\geq 0$.},colbacktitle=white,coltitle=black,boxsep=0.2mm]    
   \begin{subequations}
\footnotesize
\begin{align*}
\\[-25pt]\widecheck{s}_{lm}&\leq \underline{V}_m \underline{s} V_l+ \overline{V}_l \underline{s} V_m +\overline{V}_l \overline{V}_m\widecheck{S}-\overline{V}_l \overline{V}_m \underline{s}- \overline{V}_l \underline{V}_m \underline{s},\\
\widecheck{s}_{lm}&\leq \overline{V}_m \underline{s} V_l+\underline{V}_l \underline{s} V_m+ \overline{V}_l \overline{V}_m\widecheck{S}-\overline{V}_l \overline{V}_m \underline{s}-\underline{V}_l \overline{V}_m \underline{s},\\
\widecheck{s}_{lm}&\leq \underline{V}_m \underline{s}  V_l +  \overline{V}_l  \overline{s}  V_m +  \overline{V}_l \underline{V}_m \widecheck{S} -\overline{V}_l \underline{V}_m  \overline{s} -  \overline{V}_l \underline{V}_m  \underline{s},\\
\widecheck{s}_{lm}&\leq \overline{V}_m  \overline{s}  V_l +  \underline{V}_l  \underline{s}  V_m +  \underline{V}_l \overline{V}_m \widecheck{S} -\underline{V}_l \overline{V}_m  \overline{s} -  \underline{V}_l \overline{V}_m  \underline{s},\\
\widecheck{s}_{lm}&\leq \underline{V}_m  \overline{s}  V_l +  \overline{V}_l  \overline{s}  V_m + \underline{V}_l \underline{V}_m \widecheck{ S} -\overline{V}_l \underline{V}_m  \overline{s} -  \underline{V}_l \underline{V}_m \overline{s},\\
\widecheck{s}_{lm}&\leq \overline{V}_m \overline{s} V_l +  \underline{V}_l  \overline{s}  V_m +  \underline{V}_l \underline{V}_m \widecheck{S} -\underline{V}_l \overline{V}_m  \overline{s} -  \underline{V}_l \underline{V}_m \overline{s}.
\end{align*}
\end{subequations}
\end{tcolorbox}
}\\[0.75em]
\subfloat{%
\begin{tcolorbox}[width=\columnwidth,colback={white},title={\footnotesize
Case~III:~$ \underline{s}\geq 0$,\\[-12pt]
\begin{subequations}
\begin{align*}
 \overline{V}_l \underline{V}_m \underline{s}+ \underline{V}_l  \overline{V}_m \overline{s} &\leq  \underline{V}_l \overline{V}_m \underline{s}+ \overline{V}_l \underline{V}_m \overline{s},\\
 \overline{V}_l \underline{V}_m \underline{s}+ \underline{V}_l \overline{V}_m \overline{s} &\leq  \overline{V}_l \overline{V}_m \underline{s}+ \underline{V}_l \underline{V}_m \overline{s}.
\end{align*}
\end{subequations}},colbacktitle=white,coltitle=black,boxsep=0.2mm]    
   \begin{subequations}
\footnotesize
\begin{align*}
\\[-25pt]\widecheck{s}_{lm}&\geq \underline{V}_m \underline{s} V_l+\underline{V}_l \underline{s} V_m+\underline{V}_l \underline{V}_m \widecheck{S}-2\underline{V}_l \underline{V}_m \underline{s},\\
\widecheck{s}_{lm}&\geq \overline{V}_m \overline{s} V_l+\overline{V}_l \overline{s} V_m+\overline{V}_l \overline{V}_m\widecheck{S}-2\overline{V}_l \overline{V}_m  \overline{s},\\
\widecheck{s}_{lm}&\geq \underline{V}_m \overline{s} V_l+\underline{V}_l \overline{s} V_m + \overline{V}_l \underline{V}_m \widecheck{S}-\underline{V}_l \underline{V}_m \overline{s} -\overline{V}_l \underline{V}_m \overline{s},\\
\widecheck{s}_{lm}&\geq \overline{V}_m \underline{s} V_l+ \overline{V}_l \underline{s} V_m + \underline{V}_l \overline{V}_m \widecheck{S}-\overline{V}_l \overline{V}_m \underline{s}- \underline{V}_l \overline{V}_m \underline{s}, \\
\widecheck{s}_{lm}&\geq\frac{\Lambda_3}{ \overline{V}_l-\underline{V}_l} V_l+ \overline{V}_l \underline{s} V_m+ \overline{V}_l \underline{V}_m\widecheck{S} -\frac{\Lambda_3 \underline{V}_l}{ \overline{V}_l-\underline{V}_l} -\overline{V}_l \overline{V}_m \underline{s} \\ & \quad -\overline{V}_l \underline{V}_m \overline{s}+ \underline{V}_l \overline{V}_m \overline{s},\\
& \nonumber\text{where~} \Lambda_3 = \overline{V}_l \overline{V}_m \underline{s}-\underline{V}_l \overline{V}_m \overline{s}-\overline{V}_l \underline{V}_m \underline{s}+\overline{V}_l \underline{V}_m \overline{s},\\
\widecheck{s}_{lm}&\geq\frac{\Gamma_3}{ \underline{V}_l-\overline{V}_l} V_l+\underline{V}_l \overline{s} V_m+\underline{V}_l \overline{V}_m\widecheck{ S} -\frac{\Gamma4 \overline{V}_l}{\underline{V}_l-\overline{V}_l}-\underline{V}_l \underline{V}_m \overline{s} \\ & \quad -\underline{V}_l \overline{V}_m \underline{s}+\overline{V}_l \underline{V}_m \underline{s},\\
 & \nonumber \text{where~}  \Gamma_3 = \underline{V}_l \underline{V}_m \overline{s}-\overline{V}_l \underline{V}_m \underline{s}- \underline{V}_l \overline{V}_m \overline{s}+\underline{V}_l \overline{V}_m \underline{s}.
\end{align*}
\end{subequations}
\end{tcolorbox}
}
\vspace*{-1.5em}
\end{figure}

\begin{figure}
\subfloat{%
\begin{tcolorbox}[width=\columnwidth,colback={white},title={\footnotesize
Case~IV:~$  \overline{s} \leq 0$,\\[-12pt]
\begin{subequations}
\footnotesize
\begin{align*}
 \underline{V}_l \underline{V}_m \underline{s}+ \overline{V}_l  \overline{V}_m \overline{s} &\geq  \overline{V}_l \underline{V}_m \underline{s}+ \underline{V}_l \overline{V}_m \overline{s},\\
 \underline{V}_l \underline{V}_m \underline{s}+ \overline{V}_l \overline{V}_m \overline{s} &\geq  \underline{V}_l \overline{V}_m \underline{s}+ \overline{V}_l \underline{V}_m \overline{s}.
 \end{align*}
\end{subequations}},colbacktitle=white,coltitle=black,boxsep=0.2mm]    
   \begin{subequations}
\footnotesize
\begin{align*}
\\[-25pt] \widecheck{s}_{lm}&\leq  \underline{V}_m \overline{s}  V_l +  \underline{V}_l  \overline{s}  V_m +  \underline{V}_l \underline{V}_m\widecheck{S} -2 \underline{V}_l \underline{V}_m  \overline{s},\\
\widecheck{s}_{lm}&\leq \overline{V}_m \underline{s}  V_l +  \overline{V}_l  \underline{s}  V_m +  \overline{V}_l \overline{V}_m\widecheck{S} -2 \overline{V}_l \overline{V}_m  \underline{s},\\
\widecheck{s}_{lm}&\leq \underline{V}_m \underline{s}  V_l +  \overline{V}_l  \overline{s}  V_m +  \overline{V}_l \underline{V}_m \widecheck{S} -\overline{V}_l \underline{V}_m  \overline{s} -  \overline{V}_l \underline{V}_m  \underline{s},\\
\widecheck{s}_{lm}&\leq \overline{V}_m  \overline{s} V_l+  \underline{V}_l  \underline{s}  V_m +  \underline{V}_l \overline{V}_m \widecheck{S} -\underline{V}_l \overline{V}_m \overline{s}- \underline{V}_l \overline{V}_m \underline{s},\\
\widecheck{s}_{lm}&\leq  \underline{V}_m \underline{s}  V_l+ \underline{V}_l  \underline{s}  V_m + \frac{\Lambda_4}{\underline{s}- \overline{s}}\widecheck{S} - \frac{\Lambda_4  \overline{s}}{\underline{s}- \overline{s}}- \overline{V}_l \underline{V}_m \underline{s} \\& -\underline{V}_l \overline{V}_m \underline{s}+ \overline{V}_l  \overline{V}_m  \overline{s},\\
& \nonumber\text{where~} \Lambda_4  =  \overline{V}_l \underline{V}_m \underline{s}-\overline{V}_l \overline{V}_m \overline{s}-\underline{V}_l \underline{V}_m  \underline{s}+ \underline{V}_l \overline{V}_m  \underline{s},\\
\widecheck{s}_{lm}&\leq  \overline{V}_m \overline{s} V_l + \overline{V}_l \overline{s} V_m - \frac{\Gamma_4}{\overline{s}- \underline{s}} \widecheck{S} - \frac{\Gamma_4  \underline{s}}{ \overline{s}- \underline{s}}- \overline{V}_l \underline{V}_m  \overline{s} \\& - \underline{V}_l  \overline{V}_m  \overline{s} + \underline{V}_l  \underline{V}_m  \underline{s},\\
&\nonumber\text{where~} \Gamma_4 =  \overline{V}_l \underline{V}_m  \overline{s}- \underline{V}_l \underline{V}_m  \underline{s}- \overline{V}_l \overline{V}_m  \overline{s} + \underline{V}_l \overline{V}_m \overline{s}.
\end{align*}
\end{subequations}
\end{tcolorbox}
} \\[0.75em]
\subfloat{%
\begin{tcolorbox}[width=\columnwidth,colback={white},title={\footnotesize Case~V: ~$  \overline{s} \leq 0$,\\[-12pt] 
\begin{subequations}
\footnotesize
\begin{align*}
\overline{V}_l \underline{V}_m \underline{s}+ \underline{V}_l  \overline{V}_m \overline{s} &\geq  \underline{V}_l \overline{V}_m \underline{s}+ \overline{V}_l \underline{V}_m \overline{s},\\
\underline{V}_l \underline{V}_m \underline{s}+ \overline{V}_l  \overline{V}_m \overline{s} & <  \overline{V}_l \underline{V}_m \underline{s}+ \underline{V}_l \overline{V}_m \overline{s},\\
 \underline{V}_l \underline{V}_m \underline{s}+ \overline{V}_l \overline{V}_m \overline{s} & < \underline{V}_l \overline{V}_m \underline{s}+ \overline{V}_l \underline{V}_m \overline{s}.
\end{align*}
\end{subequations}},colbacktitle=white,coltitle=black,boxsep=0.2mm]    
   \begin{subequations}
\footnotesize
\begin{align*}
\\[-25pt]\widecheck{s}_{lm}&\leq  \underline{V}_m \overline{s}  V_l +  \underline{V}_l  \overline{s}  V_m +  \underline{V}_l \underline{V}_m \widecheck{S} -2 \underline{V}_l \underline{V}_m  \overline{s},\\
\widecheck{s}_{lm}&\leq \overline{V}_m \underline{s} V_l +  \overline{V}_l \underline{s}  V_m + \overline{V}_l \overline{V}_m\widecheck{S} -2 \overline{V}_l \overline{V}_m  \underline{s},\\
\widecheck{s}_{lm}&\leq \underline{V}_m  \underline{s}  V_l +  \underline{V}_l  \underline{s}  V_m +  \underline{V}_l \overline{V}_m \widecheck{S} -\underline{V}_l \underline{V}_m  \underline{s} -  \underline{V}_l \overline{V}_m  \underline{s},\\
\widecheck{s}_{lm}&\leq  \overline{V}_m  \overline{s}  V_l+  \overline{V}_l  \overline{s}  V_m +  \overline{V}_l \underline{V}_m \widecheck{S} -\overline{V}_l \underline{V}_m  \overline{s}-  \overline{V}_l \overline{V}_m  \overline{s},\\
\widecheck{s}_{lm}&\leq  \underline{V}_m \underline{s}  V_l + \frac{\Lambda_5}{ \underline{V}_m- \overline{V}_m}  V_m+  \overline{V}_l \underline{V}_m \widecheck{S} -\dfrac{\Lambda_5  \overline{V}_m}{\underline{V}_m- \overline{V}_m}\\&- \underline{V}_l \underline{V}_m \underline{s}  - \overline{V}_l \underline{V}_m \overline{s} + \underline{V}_l \overline{V}_m \overline{s},\\
&\nonumber\text{where~}  \Lambda_5 =  \underline{V}_l \underline{V}_m  \underline{s}-\underline{V}_l \overline{V}_m \overline{s}-\overline{V}_l \underline{V}_m  \underline{s}+ \overline{V}_l \underline{V}_m  \overline{s},\\
\widecheck{s}_{lm}&\leq  \overline{V}_m \overline{s}  V_l + \frac{\Gamma_5}{ \overline{V}_m- \underline{V}_m}  V_m+  \underline{V}_l \overline{V}_m \widecheck{S} -\frac{\Gamma_5  \underline{V}_m}{ \overline{V}_m- \underline{V}_m}\\&- \underline{V}_l \overline{V}_m  \underline{s}  -  \overline{V}_l  \overline{V}_m  \overline{s} +  \overline{V}_l  \underline{V}_m  \underline{s},\\
&\nonumber \text{where~}\Gamma_5 =  \underline{V}_l \overline{V}_m  \underline{s}- \overline{V}_l \underline{V}_m  \underline{s}- \underline{V}_l \overline{V}_m \overline{s} + \overline{V}_l \overline{V}_m  \overline{s}.
\end{align*}
\end{subequations}
\end{tcolorbox}
}
\end{figure}

\begin{figure}
\subfloat{%
\begin{tcolorbox}[width=\columnwidth,colback={white},title={\footnotesize Case VI:~$ \underline{s} \leq 0, \overline{s} \geq 0$.},colbacktitle=white,coltitle=black,boxsep=0.2mm]    
   \begin{subequations}
\footnotesize
\begin{align*}
\\[-25pt]\widecheck{s}_{lm}& \geq  \overline{V}_m \overline{s}  V_l +  \overline{V}_l  \overline{s}  V_m +  \overline{V}_l \overline{V}_m \widecheck{S} -2 \overline{V}_l \overline{V}_m  \overline{s},\\
\widecheck{s}_{lm}&\geq \overline{V}_m  \underline{s}  V_l +  \underline{V}_l  \overline{s}  V_m +  \underline{V}_l \overline{V}_m \widecheck{S} -\underline{V}_l \overline{V}_m \underline{s} -  \underline{V}_l \overline{V}_m  \overline{s},\\
\widecheck{s}_{lm}&\geq \overline{V}_m  \underline{s}  V_l +  \underline{V}_l  \underline{s}  V_m +  \underline{V}_l \underline{V}_m \widecheck{S} -\underline{V}_l \overline{V}_m  \underline{s} -  \underline{V}_l \underline{V}_m  \underline{s},\\
\widecheck{s}_{lm}&\geq  \underline{V}_m  \overline{s}  V_l+  \overline{V}_l  \underline{s}  V_m +  \overline{V}_l \underline{V}_m \widecheck{S} -\overline{V}_l \underline{V}_m  \overline{s}-  \overline{V}_l \underline{V}_m  \underline{s},\\
\widecheck{s}_{lm}&\geq  \underline{V}_m  \underline{s}  V_l+  \overline{V}_l  \underline{s}  V_m +  \underline{V}_l \underline{V}_m \widecheck{S} -\overline{V}_l \underline{V}_m  \underline{s}-   \underline{V}_l \underline{V}_m  \underline{s},\\
\widecheck{s}_{lm}&\geq  \underline{V}_m \overline{s}  V_l +  \underline{V}_l  \overline{s}  V_m + \frac{\Lambda_6}{\overline{s}- \underline{s}}  \widecheck{S}-\frac{\Lambda_6  \underline{s}}{ \overline{s}- \underline{s}}- \underline{V}_l \overline{V}_m  \overline{s} \\& -  \overline{V}_l  \underline{V}_m  \overline{s} +  \overline{V}_l  \overline{V}_m  \underline{s},\\
& \nonumber\text{where~}  \Lambda_6 =  \underline{V}_l \overline{V}_m  \overline{s}-\overline{V}_l \overline{V}_m  \underline{s}-\underline{V}_l \underline{V}_m  \overline{s}+ \overline{V}_l \underline{V}_m  \overline{s}.
\end{align*}
\end{subequations}
\end{tcolorbox}
}\\[0.75em]
\subfloat{%
\begin{tcolorbox}[width=\columnwidth,colback={white},title={\footnotesize Case~VII:~$ \underline{s} \leq 0,  \overline{s} \geq 0$.},colbacktitle=white,coltitle=black,boxsep=0.2mm]    
   \begin{subequations}
\footnotesize
\begin{align*}
\\[-25pt] \widecheck{s}_{lm}&\leq  \overline{V}_m \underline{s}  V_l +  \overline{V}_l  \underline{s}  V_m +  \overline{V}_l \overline{V}_m \widecheck{ S} -2 \overline{V}_l \overline{V}_m  \underline{s},\\
\widecheck{s}_{lm} &\leq \underline{V}_m  \underline{s}  V_l +  \overline{V}_l \overline{s} V_m + \overline{V}_l \underline{V}_m \widecheck{S} -\overline{V}_l \underline{V}_m \overline{s} - \overline{V}_l \underline{V}_m  \underline{s},\\
\widecheck{s}_{lm} &\leq \overline{V}_m  \overline{s}  V_l +  \underline{V}_l  \overline{s}  V_m +  \underline{V}_l \underline{V}_m \widecheck{S} -\underline{V}_l \overline{V}_m  \overline{s} -  \underline{V}_l \underline{V}_m  \overline{s},\\
\widecheck{s}_{lm} &\leq \overline{V}_m \overline{s} V_l+  \underline{V}_l  \underline{s}  V_m +  \underline{V}_l \overline{V}_m \widecheck{S} -\underline{V}_l \overline{V}_m  \overline{s}-  \underline{V}_l \overline{V}_m  \underline{s},\\
\widecheck{s}_{lm} &\leq  \underline{V}_m  \overline{s}  V_l+  \overline{V}_l  \overline{s}  V_m +  \underline{V}_l \underline{V}_m \widecheck{S} -\overline{V}_l \underline{V}_m  \overline{s} - \underline{V}_l \underline{V}_m  \overline{s},\\
\widecheck{s}_{lm} &\leq  \underline{V}_m \underline{s}  V_l +  \underline{V}_l  \underline{s}  V_m + \frac{\Lambda_7}{ \underline{s}- \overline{s}}  \widecheck{S}-\frac{\Lambda_7 \overline{s}}{ \underline{s}- \overline{s}} - \overline{V}_l \underline{V}_m  \underline{s} \\ & \qquad - \underline{V}_l  \overline{V}_m  \underline{s} +  \overline{V}_l  \overline{V}_m  \overline{s},\\
 & \nonumber\text{where~}\Lambda_7 =  \overline{V}_l \underline{V}_m  \underline{s}-\overline{V}_l \overline{V}_m  \overline{s}-\underline{V}_l \underline{V}_m  \underline{s}+ \underline{V}_l \overline{V}_m  \underline{s}.
\end{align*}
\end{subequations}
\end{tcolorbox}
}
\vspace*{-1.5em}
\end{figure}


\bibliographystyle{IEEEtran}
\bibliography{ref}

\begin{thebibliography}{10}
\providecommand{\url}[1]{#1}
\csname url@samestyle\endcsname
\providecommand{\newblock}{\relax}
\providecommand{\bibinfo}[2]{#2}
\providecommand{\BIBentrySTDinterwordspacing}{\spaceskip=0pt\relax}
\providecommand{\BIBentryALTinterwordstretchfactor}{4}
\providecommand{\BIBentryALTinterwordspacing}{\spaceskip=\fontdimen2\font plus
\BIBentryALTinterwordstretchfactor\fontdimen3\font minus
  \fontdimen4\font\relax}
\providecommand{\BIBforeignlanguage}[2]{{%
\expandafter\ifx\csname l@#1\endcsname\relax
\typeout{** WARNING: IEEEtran.bst: No hyphenation pattern has been}%
\typeout{** loaded for the language `#1'. Using the pattern for}%
\typeout{** the default language instead.}%
\else
\language=\csname l@#1\endcsname
\fi
#2}}
\providecommand{\BIBdecl}{\relax}
\BIBdecl

\bibitem{bukhsh_tps}
W.~Bukhsh, A.~Grothey, K.~McKinnon, and P.~Trodden, ``{Local Solutions of the
  Optimal Power Flow Problem},'' \emph{IEEE Trans. Power Syst.}, vol.~28,
  no.~4, pp. 4780--4788, 2013.

\bibitem{pascalNPhard}
K.~Lehmann, A.~Grastien, and P.~{Van Hentenryck}, ``{AC-Feasibility on Tree
  Networks is NP-Hard},'' \emph{IEEE Trans. Power Syst.}, vol.~31, no.~1, pp.
  798--801, January 2016.

\bibitem{bienstock2015nphard}
D.~Bienstock and A.~Verma, ``{Strong NP-hardness of AC Power Flows
  Feasibility},'' \emph{arXiv:1512.07315}, Dec. 2015.

\bibitem{opf_litreview1993IandII}
J.~Momoh, R.~Adapa, and M.~El-Hawary, ``{A Review of Selected Optimal Power
  Flow Literature to 1993. Parts I and II},'' \emph{IEEE Trans. Power Syst.},
  vol.~14, no.~1, pp. 96--111, Feb. 1999.

\bibitem{ferc4}
A.~Castillo and R.~O'Neill, ``{Survey of Approaches to Solving the ACOPF (OPF
  Paper 4)},'' FERC, Tech. Rep., Mar. 2013.

\bibitem{marley2016}
J.~F. Marley, D.~K. Molzahn, and I.~A. Hiskens, ``{Solving Multiperiod OPF
  Problems using an AC-QP Algorithm Initialized with an SOCP Relaxation},''
  \emph{IEEE Trans. Power Syst.}, vol.~32, no.~5, pp. 3538--3548, Sept. 2017.

\bibitem{molzahn2017survey}
D.~K. Molzahn and I.~A. Hiskens, ``{A Survey of Relaxations and Approximations
  of the Power Flow Equations},'' \emph{{\rm invited submission to} Found.
  Trends Electric Energy Syst.}, 2018.

\bibitem{coffrin2015qc}
C.~Coffrin, H.~Hijazi, and P.~{Van Hentenryck}, ``{The QC Relaxation: A
  Theoretical and Computational Study on Optimal Power Flow},'' \emph{IEEE
  Trans. Power Syst.}, vol.~31, no.~4, pp. 3008--3018, July 2016.

\bibitem{coffrin2016strengthen_tps}
C.~Coffrin, H.~L. Hijazi, and P.~{Van Hentenryck}, ``{Strengthening the SDP
  Relaxation of AC Power Flows with Convex Envelopes, Bound Tightening, and
  Valid Inequalities},'' \emph{IEEE Trans. Power Syst.}, vol.~32, no.~5, pp.
  3549--3558, Sept. 2017.

\bibitem{chen2015}
C.~Chen, A.~Atamt{\"u}rk, and S.~Oren, ``{Bound Tightening for the Alternating
  Current Optimal Power Flow Problem},'' \emph{IEEE Trans. Power Syst.},
  vol.~31, no.~5, pp. 3729--3736, Sept. 2016.

\bibitem{StrongSOCPRelaxations}
B.~Kocuk, S.~S. Dey, and X.~A. Sun, ``{Strong SOCP Relaxations for the Optimal
  Power Flow Problem},'' \emph{Oper. Res.}, vol.~64, no.~6, pp. 1177--1196, May
  2016.

\bibitem{arctan2}
------, ``{Matrix Minor Reformulation and SOCP-based Spatial Branch-and-Cut
  Method for the AC Optimal Power Flow Problem},'' \emph{arXiv:1703.03050},
  March 2017.

\bibitem{chen2015cuts}
C.~Chen, A.~Atamt{\"u}rk, and S.~S. Oren, ``{A Spatial Branch-and-Cut Algorithm
  for Nonconvex QCQP with Bounded Complex Variables},'' \emph{Math. Prog.}, pp.
  1--29, 2016.

\bibitem{coffrin2016quadtrig}
K.~Bestuzheva, H.~L. Hijazi, and C.~Coffrin, ``{Convex Relaxations for
  Quadratic On/Off Constraints and Applications to Optimal Transmission
  Switching},'' \emph{{\rm Preprint:}
  \url{http://www.optimization-online.org/DB_FILE/2016/07/5565.pdf}}, 2016.

\bibitem{ruiz2011}
J.~Ruiz and J.~Grossmann, ``{Using Redundancy to Strengthen the Relaxation for
  the Global Optimization of MINLP Problems},'' \emph{Comput. \& Chemical
  Eng.}, vol.~35, no.~12, 2011.

\bibitem{mccormick1976}
G.~{McCormick}, ``{Computability of Global Solutions to Factorable Nonconvex
  Programs: Part I--Convex Underestimating Problems},'' \emph{Math. Prog.},
  vol.~10, no.~1, pp. 147--175, 1976.

\bibitem{meyer2004a}
C.~Meyer and C.~Floudas, \emph{{Trilinear Monomials with Positive or Negative
  Domains: Facets of the Convex and Concave Envelopes}}.\hskip 1em plus 0.5em
  minus 0.4em\relax Boston, MA: Springer US, 2004, pp. 327--352.

\bibitem{meyer2004b}
------, ``{Trilinear Monomials with Mixed Sign Domains: Facets of the Convex
  and Concave Envelopes},'' \emph{J. Global Optimiz.}, vol.~29, no.~2, pp.
  125--155, 2004.

\bibitem{matpower}
R.~Zimmerman, C.~Murillo-S{\'a}nchez, and R.~Thomas, ``{MATPOWER: Steady-State
  Operations, Planning, and Analysis Tools for Power Systems Research and
  Education},'' \emph{IEEE Trans. Power Syst.}, no.~99, pp. 1--8, 2011.

\bibitem{Jabr2006}
R.~Jabr, ``{Radial Distribution Load Flow using Conic Programming},''
  \emph{IEEE Trans. Power Syst.}, vol.~21, no.~3, pp. 1458--1459, Aug. 2006.

\bibitem{molzahn-opf_spaces}
D.~K. Molzahn, ``{Computing the Feasible Spaces of Optimal Power Flow
  Problems},'' \emph{IEEE Trans. Power Syst.}, vol.~32, no.~6, pp. 4752--4763,
  Nov. 2017.

\bibitem{nesta3}
C.~Coffrin, D.~Gordon, and P.~Scott, ``{NESTA, the NICTA Energy System Test
  Case Archive (v0.7)},'' \emph{arXiv:1411.0359}, June 2017.

\bibitem{taylorlinearRelaxationOPF}
J.~Taylor and F.~Hover, ``{Linear Relaxations for Transmission System
  Planning},'' \emph{IEEE Trans. Power Syst.}, vol.~26, no.~4, pp. 2533--2538,
  Nov. 2011.

\bibitem{ACCpaper}
M.~R. Narimani, D.~K. Molzahn, D.~Wu, and {M. L. Crow}, ``{Empirical
  Investigation of Non-Convexities in Optimal Power Flow Problems},''
  \emph{{\rm to appear in} American Control Conf. (ACC)}, June 2018.

\bibitem{yalmip}
J.~Lofberg, ``{YALMIP: A Toolbox for Modeling and Optimization in MATLAB},'' in
  \emph{{IEEE Int. Symp. Compu. Aided Control Syst. Des.}}, 2004, pp. 284--289.

\bibitem{kocuk2017_ots}
B.~Kocuk, S.~S. Dey, and X.~A. Sun, ``{New Formulation and Strong MISOCP
  Relaxations for AC Optimal Transmission Switching Problem},'' \emph{IEEE
  Trans. Power Syst.}, vol.~32, no.~6, pp. 4161--4170, Nov. 2017.

\bibitem{molzahn2018pscc}
D.~K. Molzahn and L.~A. Roald, ``{Towards an AC Optimal Power Flow Algorithm
  with Robust Feasibility Guarantees},'' \emph{20th Power Syst. Comput. Conf.
  (PSCC)}, June 2018.

\end{thebibliography}

 \begin{landscape}
 \begin{table}[t]
\caption{Results from Applying the QC Relaxation with Various Improvements to Selected Test Cases}
\hspace*{2.5em}\resizebox{6.5in}{!}{\begin{minipage}{\textwidth}

 \newpage 
 
\footnotesize

\label{tab:QC_results}%

\begin{tabular}{|l||c||c|c|c|c|c|c||c|c|c|c|c|c|c|c|c|c|}
\hline 
 &  & \multicolumn{6}{c||}{\thead{Optimality Gap (\%)}} & \multicolumn{2}{c|}{\thead{All \\Cons.}} & \multicolumn{2}{c|}{\thead{w/o\\ MF}} & \multicolumn{2}{c|}{\thead{w/o\\ $\bm{\Delta}$}} & \multicolumn{2}{c|}{\thead{w/o \\ \{MF, $\bm{\Delta}$\}}} & \thead{w/o\\ BT} & \thead{w/o \{MF, \\ BT, $\bm{\Delta}$\}}\tabularnewline
\hline 
\thead{Test Case} & \thead{AC \\(\$/hr)} & \thead{All \\ Cons.} & \thead{w/o \\MF} & \thead{w/o\\ $\bm{\Delta}$} & \thead{w/o \\ \{MF, $\bm{\Delta}$\}} & \thead{w/o BT} & \thead{w/o \{BT, \\ MF, $\bm{\Delta}$\}} & \thead{BT\\ time} & \thead{QC\\ time} & \thead{BT \\time} & \thead{QC \\time} & \thead{BT \\time} & \thead{QC \\time} & \thead{BT \\time}& \thead{QC \\time} & \thead{QC \\time} & \thead{QC \\time}\tabularnewline
\hline 
\hline 
nesta\_case3\_lmbd & \hphantom{00}5812.64 & \hphantom{0}0.17 & \hphantom{0}0.21 & \hphantom{0}0.17 & \hphantom{0}0.21 & \hphantom{0}1.13 & \hphantom{0}1.23 & \hphantom{000}4.0 & 0.4 &\hphantom{000}4.1 & 0.4 & \hphantom{00}3.5 & 0.4 & \hphantom{0}11.5 & 0.6 & \hphantom{0}0.4 & \hphantom{0}0.2\tabularnewline
\hline 
nesta\_case5\_pjm &  \hphantom{0}17551.89 & 11.57 & 15.39 & 11.63 & 15.50 & 17.00 & 17.01 & \hphantom{000}6.1 & 0.5 & \hphantom{000}5.9 & 0.4 & \hphantom{00}4.9 & 0.4 & \hphantom{0}10.4 & 0.5 & \hphantom{0}0.3 & \hphantom{0}0.2\tabularnewline
\hline 
nesta\_case29\_edin &  \hphantom{0}29895.49 & \hphantom{0}0.03 & \hphantom{0}0.05 & \hphantom{0}0.02 & \hphantom{0}0.05 & \hphantom{0}0.09 & \hphantom{0}0.10 & \hphantom{0}170.8 & 0.6 & \hphantom{0}175.6 & 0.4 & 104.8 & 0.5 & 137.8 & 0.4 & \hphantom{0}0.6 & \hphantom{0}0.3\tabularnewline
\hline 
nesta\_case118\_ieee &  \hphantom{00}3718.64 & \hphantom{0}0.35 & \hphantom{0}0.41 & \hphantom{0}0.36 & \hphantom{0}0.42 & \hphantom{0}1.38 & \hphantom{0}1.64 & \hphantom{0}943.2 & 1.1 & \hphantom{0}799.4 & 0.9 & 756.4 & 1.0 & 648.9 & 0.9 & \hphantom{0}2.1 & \hphantom{0}1.5\tabularnewline
\hline 
nesta\_case29\_edin\_\_api & \hphantom{}295291.22 & \hphantom{0}0.07 & \hphantom{0}0.09 & \hphantom{0}0.08 & \hphantom{0}0.09 & \hphantom{0}0.40 & \hphantom{0}0.40 & \hphantom{0}133.9 & 0.5 & \hphantom{0}126.4 & 0.5 & 167.1 & 0.5 & 106.9 & 0.7 & \hphantom{0}0.5 & \hphantom{0}0.4\tabularnewline
\hline 
nesta\_case30\_fsr\_\_api & \hphantom{000}366.57 & \hphantom{0}4.01 & \hphantom{0}5.29 & \hphantom{0}4.16 & \hphantom{0}5.50 & 82.47 & 82.46 & \hphantom{0}102.1 & 0.6 & \hphantom{00}74.9 & 0.5 & \hphantom{0}72.3 & 0.5 & \hphantom{0}59.6 & 0.8 & \hphantom{0}0.4 & \hphantom{0}0.3\tabularnewline
\hline 
nesta\_case73\_ieee\_rts\_\_api & \hphantom{0}19995.00 & \hphantom{0}0.11 & \hphantom{0}0.18 & \hphantom{0}0.11 & \hphantom{0}0.19 & 10.11 & 16.50 & \hphantom{0}317.2 & 0.7 & \hphantom{0}326.6 & 0.7 & 197.3 & 0.7 & 220.8 & 0.8 & \hphantom{0}0.7 & \hphantom{0}0.4\tabularnewline
\hline 
nesta\_case118\_ieee\_\_api & \hphantom{0}10269.82 & 13.61 & 18.32 & 13.90 & 18.71 & 76.47 & 77.06 & 1757.8 & 0.9 & 1439.0 & 0.6 & 726.6 & 0.7 & 603.6 & 0.5 & \hphantom{0}0.9 & \hphantom{0}0.6\tabularnewline
\hline 
nesta\_case3\_lmbd\_\_sad & \hphantom{00}5959.33 & \hphantom{0}0.03 & \hphantom{0}0.03 & \hphantom{0}0.11 & \hphantom{0}0.11 & \hphantom{0}1.31 & \hphantom{0}1.34 & \hphantom{000}5.8 & 0.4 & \hphantom{000}5.2 & 0.4 & \hphantom{00}5.2 & 0.5 & \hphantom{0}12.1 & 0.4 & \hphantom{0}0.3 & \hphantom{0}0.2\tabularnewline
\hline 
nesta\_case5\_pjm\_\_sad & \hphantom{0}26115.20 & \hphantom{0}0.10 & \hphantom{0}0.12 & \hphantom{0}0.10 & \hphantom{0}0.12 & \hphantom{0}0.78 & \hphantom{0}1.35 & \hphantom{000}5.5 & 0.5 & \hphantom{000}4.7 & 0.4 & \hphantom{00}4.8 & 0.5 & \hphantom{00}5.7 & 0.4 & \hphantom{0}0.3 & \hphantom{0}0.2\tabularnewline
\hline 
nesta\_case9\_wscc\_\_sad & \hphantom{00}5528.26 & \hphantom{0}0.04 & \hphantom{0}0.05 & \hphantom{0}0.04 & \hphantom{0}0.05 & \hphantom{0}0.47 & \hphantom{0}0.54 & \hphantom{000}4.6 & 0.5 & \hphantom{000}4.0 & 0.4 & \hphantom{00}3.9 & 0.4 & \hphantom{00}4.4 & 0.4 & \hphantom{0}0.4 & \hphantom{0}0.2\tabularnewline
\hline 
nesta\_case24\_ieee\_rts\_\_sad & \hphantom{0}76943.25 & \hphantom{0}0.07 & \hphantom{0}0.09 & \hphantom{0}0.07 & \hphantom{0}0.09 & \hphantom{0}2.74 & \hphantom{0}3.15 & \hphantom{00}54.9 & 0.6 & \hphantom{00}45.2 & 0.4 & \hphantom{0}38.7 & 0.5 & \hphantom{0}76.5 & 0.5 & \hphantom{0}0.5 & \hphantom{0}0.3\tabularnewline
\hline 
nesta\_case29\_edin\_\_sad & \hphantom{0}41258.49 & \hphantom{0}1.75 & \hphantom{0}2.72 & \hphantom{0}1.79 & \hphantom{0}2.68 & 19.64 & 27.94 & \hphantom{0}220.5 & 0.6 & \hphantom{0}631.9 & 0.5 & 133.4 & 0.5 & 357.2 & 0.7 & \hphantom{0}0.5 & \hphantom{0}0.3\tabularnewline
\hline 
nesta\_case30\_as\_\_sad & \hphantom{000}897.49 & \hphantom{0}0.14 & \hphantom{0}0.17 & \hphantom{0}0.14 & \hphantom{0}0.17 & \hphantom{0}2.30 & \hphantom{0}2.37 & \hphantom{00}39.4 & 0.5 & \hphantom{00}30.0 & 0.4 & \hphantom{0}27.0 & 0.5 & \hphantom{0}47.2 & 0.7 & \hphantom{0}0.4 & \hphantom{0}0.3\tabularnewline
\hline 
nesta\_case39\_epri\_\_sad & \hphantom{0}96745.01 & \hphantom{0}0.00 & \hphantom{0}0.01 & \hphantom{0}0.01 & \hphantom{0}0.02 & \hphantom{0}0.05 & \hphantom{0}0.08 & \hphantom{0}101.0 & 0.6 & \hphantom{00}75.7 & 0.5 & \hphantom{0}66.5 & 0.5 & \hphantom{0}57.7 & 0.9 & \hphantom{0}0.5 & \hphantom{0}0.3\tabularnewline
\hline 
nesta\_case118\_ieee\_\_sad & \hphantom{00}4106.72 & \hphantom{0}0.67 & \hphantom{0}0.90 & \hphantom{0}0.68 & \hphantom{0}0.92 & \hphantom{0}4.10 & \hphantom{0}4.96 & \hphantom{0}937.6 & 0.8 & \hphantom{0}906.7 & 0.6 & 996.7 & 0.7 & 733.2 & 1.2 & \hphantom{0}0.6 & \hphantom{0}0.6\tabularnewline
\hline 
nmwc14 & \hphantom{00}2529.87 & \hphantom{0}0.17 & \hphantom{0}0.17 & \hphantom{0}0.19 & \hphantom{0}0.19 & \hphantom{0}0.22 & \hphantom{0}0.22 & \hphantom{00}26.7 & 1.5 & \hphantom{00}21.9 & 1.3 & \hphantom{0}16.8 & 1.4 & \hphantom{0}14.9 & 1.3 & \hphantom{0}2.1 & \hphantom{0}2.0\tabularnewline
\hline 
nmwc57 & \hphantom{00}9186.12 & \hphantom{0}6.44 & \hphantom{0}6.44 & \hphantom{0}7.22 & \hphantom{0}7.20 & \hphantom{0}9.66 & \hphantom{0}9.66 & \hphantom{0}322.1 & 2.0 & \hphantom{0}266.2 & 1.8 & 352.0 & 2.1 & 232.1 & 1.9 & 11.9 & 11.7\tabularnewline
\hline 
nmwc118 & \hphantom{0}34663.69 & 17.07 & 17.26 & 23.15 & 23.19 & 24.07 & 24.07 & 2124.4 & 3.0 & 1624.1 & 2.7 & 928.9 & 2.9 & 542.4 & 2.6 & \hphantom{0}3.3 & \hphantom{0}3.0\tabularnewline
\hline
\multicolumn{18}{l}{Abbreviations: ``All Constraints'' (\emph{All Cons.}), ``without'' (\emph{w/o}), ``Meyer and Floudas Envelopes'' (\emph{MF}), ``Voltage Magnitude Difference constraints'' ($\Delta$), and ``Bound Tightening'' (\emph{BT}). All times in seconds.} \tabularnewline
\end{tabular}

 \end{minipage}}
 \end{table}
\end{landscape}

\end{document}